\documentclass{article}
\usepackage{mathptmx}       % selects Times Roman as basic font
\usepackage{helvet}         % selects Helvetica as sans-serif font
\usepackage{courier}        % selects Courier as typewriter font
\usepackage{type1cm}        % activate if the above 3 fonts are not available on your system
\usepackage{makeidx}         % allows index generation
\usepackage{graphicx}        % standard LaTeX graphics tool when including figure files
\usepackage{multicol}        % used for the two-column index
\usepackage[bottom]{footmisc}% places footnotes at page bottom
\usepackage{subfigure}
\usepackage{amsfonts}
\usepackage{amssymb}
\usepackage[cmex10]{amsmath}
\usepackage{epstopdf}
\usepackage{authblk}
\usepackage{enumitem} 
\usepackage{sectsty}

\sectionfont{\large}
\subsectionfont{\small}
\date{}

\begin{document}
	
	\title{Eigenvalues of random networks with cycles}
	
	\author{Pau Vilimelis Aceituno }
	\affil{Max Planck Institute for Mathematics in the Sciences\\
		Max Planck School of Cognition,\\
		Leipzig, Germany
		\\ aceituno@mis.mpg.de
	}

	\maketitle    

\begin{abstract}
Networks are often studied using the eigenvalues of their adjacency matrix, a powerful mathematical tool with a wide range of applications. Since in real systems the exact graph structure is not known, researchers resort to random graphs to obtain eigenvalue properties from known structural features. However, this theory is far from intuitive and often requires training of free probability, cavity methods or a strong familiarity with probability theory. In this note we offer a different perspective on this field by focusing on the cycles in a graph. We use the so-called method of moments to obtain relation between eigenvalues and cycle weights and then we obtain spectral properties of random graphs with cyclic motifs. We use it to explore properties of the eigenvalues of adjacency matrices of graphs with short cycles and of circulant directed graphs. Although our result is not as powerful as the some of the existing methods, they are nevertheless useful and far easier to understand.
\end{abstract}

\section{Introduction}
Networks and graphs are a very flexible modeling tool that has been crucial in studying complex systems. Of particular importance in those studies are the spectral methods, a family of tools that use the eigenvalues of the adjacency matrix of the graph \cite{spielman2007spectral, newman2010networks} and which has been applied to fields as diverse as wireless communications \cite{tulino2004random}, ecology \cite{allesina2012stability}, biological networks \cite{luo2006application}, image processing \cite{shi2000normalized}, machine learning  \cite{aceituno2017tailoring} or nuclear physics \cite{haq1982fluctuation}. However, in real world systems with many interacting elements the exact graph is not known, and scientists must estimate those eigenvalues in order to use spectral methods. Hence researchers resort to spectral random matrix theory, the mathematical field that studies the distribution of eigenvalues from a probabilistic perspective \cite{edelman2005random}. However, this approach requires a deep understanding of probability theory \cite{tao2012topics} or of statistical mechanics \cite{spencer2012duality}, and its results can be difficult to understand for non-experts. In our work is based so-called method of moments used to prove Winger's semicircle law \cite{wigner1993characteristic}, that has been extended to other families of symmetric matrices in an intuitive way \cite{kirsch2016sixty}. However, we give an extra twist by giving it a graph-theoretical interpretation, allowing us to go beyond the symmetric case and providing an intuitive approach to network spectra.

The rest of this paper is organized as follows: First we explain the connection between cycles in a graph and the eigenvalues of its adjacency matrix, then we use this connection to study the symmetries embedded in the spectra of random graphs with circular motifs and finally the spectral radius of directed circulant graphs.

\section{Framework}
A directed graph $G(N,E)$ is a set of $N$ nodes and $E$ edges, where every edge $e=(n,m)$ represents a directed connection going from node $n$ to node $m$. Each edge has an associated weight $w(e)$ which represents the strength of that connection and might be negative. The graph can be represented by its adjacency matrix $M$, an $N\times N$ matrix where every entry $M_{nm}$ is the weight of edge $(n,m)$, and non-existent edges correspond to entries with value zero.

We define the normalized weight of the cycles of length $L$
\begin{equation}\label{eq:defcL}
\rho_L = \dfrac{1}{N}\sum_{c\in C_L}w_c
\end{equation}
where $C_L$ is the set of cycles of length $L$, and $w_c = \prod_{e\in c} w(e)$, the multiplication of weights of the edges $e$ in cycle $c$. Note that this includes cycles where edges or nodes are visited multiple times, as opposed to simple cycles where every node and edge can only be visited once.
Since the value $\sum_{c\in C_L}w_c$ is given by the entries of the power of the graphs' adjacency matrix $M$ \cite{gross2005graph}, we can obtain $\rho_L$ from the adjacency matrix $M$ by the formula
\begin{equation}\label{eq:rhoTrace}
\rho_L = \dfrac{1}{N}\text{tr}\left[M^L\right]
\end{equation}
where $\text{tr}\left[\cdot\right]$ is the trace operator. The trace of a matrix equals the sum of its eigenvalues \cite{strang1993introduction}, therefore $\rho_L$ can be written as
\begin{equation}\label{eq:cLIsSumOfEigs}
\rho_L = \dfrac{1}{N}\sum_{n=1}^{N} \lambda_n^L
\end{equation} 
where $\lambda_n$ is the $n$th eigenvalue of the adjacency matrix $M$. 

Equation \ref{eq:cLIsSumOfEigs} is particularly interesting as $N\to \infty$. For this limit we can consider the eigenvalues of a matrix as random i.i.d. values sampled from a probability density function $p(\lambda)$ in the complex plane $\mathbb{C}$. For a random graph sampled from a probability distribution with a fixed $\rho_L$, this value corresponds to a moment of the eigenvalue distribution,
\begin{equation}\label{eq:LthMoment}
	\mu_L = \int_\mathbb{C} p(\lambda) \lambda^L d\lambda = \lim\limits_{N\to \infty} \dfrac{1}{N}\sum_{n=1}^{N} \lambda_n^L = \rho_L.
\end{equation}

\section{Random Networks with Cyclic Motifs}

Many real-world networks have an unexpectedly high number of network motifs, small subgraphs that are overrepresented with respect to randomized networks \cite{milo2002network}. One of the most relevant examples is the feedback loop, a well-known structure in control theory \cite{doyle2013feedback} which is equivalent to a directed cycle in a network. This section shows how feedback loops shape graph spectra by studying graphs with abundant cyclic motifs.
We will consider sparse random graphs sampled from a distribution where all the connections between nodes have the same probability and with weights independently drawn from $\left[-w, w\right]$ with equal probability. In order to work on the limit $N\rightarrow \infty$, we will normalize the weights so that the variance of incoming and outgoing edge weights is one, therefore $w = \sqrt{\frac{N}{E}}$. To account for the feedback motifs, a fraction of the edges are embedded into cycles of length $\tau$ where the multiplication of the weights of edges in each cycle is either always positive or always negative, corresponding, respectively to positive or negative feedback. 

As the number of nodes goes to infinity, the values of $\rho_L$ converge to their expectations, and since only the cycles of length $\tau$ have non-zero expectations. 
\begin{equation}\label{eq:motifsRhoTau}
	\lim\limits_{N\to \infty} \rho_L =
	\begin{cases}
	w^L \dfrac{F_L}{N} \ \ \iff L \equiv 0 \mod{\tau}\\
	 0 \quad \quad \text{otherwise}
	\end{cases}
\end{equation}
Where $C_L$ is the number of cycles of length $L$ that were added as feedback loops. Since $\rho_L$ are the moments of $p(\lambda)$, the previous equation says that the moments which are not multiples of $\tau$ tend to zero. We will use this to study the symmetries in the eigenvalues, but first we need to consider the rotated density function $p(e^{\theta i}\lambda)$. The moments of this distribution can be computed through a simple change of variable,
\begin{equation}
	\mu_L^\theta = \int_{\mathbb{C}} p(e^{\theta i}\lambda) \lambda^L d\lambda 
	= \int_{\mathbb{C}}  p(\lambda) \lambda^L e^{L\theta i} d\lambda 
	= e^{L\theta i}\int_{\mathbb{C}}  p(\lambda) \lambda^L d\lambda 
	= e^{L\theta i}\mu_L.
\end{equation}
Then the moments of $p(\lambda)$ and $p(e^{\theta i}\lambda)$ are equal under the following conditions,
\begin{equation}
\mu_L = \mu_L^\theta = e^{L\theta i}\mu_L
\iff \begin{cases}
\theta \in \left\{0, \dfrac{2\pi}{L}, \dfrac{4\pi}{L},..., \dfrac{2(L-1)\pi}{L}\right\}
\\
\mu_L = 0
\end{cases}
\end{equation}
By using Eq. \ref{eq:motifsRhoTau}, in a large graph with abundant cycles of length $\tau$,
\begin{equation}
	\forall L\ \ \mu_L^\theta = \mu_L\ \ \ 
	 \iff \theta = 0, \dfrac{2\pi}{\tau}, \dfrac{4\pi}{\tau}...
\end{equation}
And since the eigenvalues are bounded, the equality of moments implies that $p(\lambda) = p(e^{\theta i}\lambda)$, so 
\begin{equation}\label{eq:DensitySymmetry}
	p(\lambda) = p(e^{\theta i}\lambda)
\end{equation}
for $\theta \in \left\{0, \frac{2\pi}{\tau},..., \frac{2(\tau-1)\pi}{\tau}\right\}$. In geometric terms, this means that $p(\lambda)$ has $\tau$ rotational symmetries in the complex plane. 

Since the entries of our adjacency matrix are real, the eigenvalue distribution is also symmetric with respect to the real line. Combining this axis of symmetry with the rotational symmetries from Eq. \ref{eq:DensitySymmetry}, we obtain that $p(\lambda)$ has $\tau$ axes of symmetry in the complex plane, all of them passing through the origin and with at angles $0, \frac{2 \pi}{\tau}, ...,  \pi-\frac{2 \pi }{\tau}$. 

Within the constraints of this symmetries, we would still need to know how $\rho_\tau$ affects the eigenvalue distribution. 

Since the entries of our adjacency matrix are i.i.d. and their magnitude decreases sub-exponentially, the probability distributions are uniform \cite{bourgade2014local} within a closed curve. This means that increasing $\rho_\tau$ decreases the surface in the angular range $\theta \in \left[ \frac{\pi}{k}, \frac{2\pi}{k}\right]$ while increasing it in $\theta \in \left[ 0, \frac{\pi}{k}\right]$. Intuitively, $\rho_\tau>0$ accounts for the number of cycles, and every cycle is a feedback loop, therefore the linear system described by $M$ has many cycles with positive feedback, so it resonates at the frequency $\frac{1}{\tau}$ and its corresponding harmonics. This resonant frequency corresponds to the phases of the dominant eigenvalues, meaning that when $\rho_\tau$ increases the eigenvalues with the phases $\frac{2\pi k}{\tau}$ for $k \in \{0,1,...,\tau-1\}$ increase their moduli while those with phases $\frac{2\pi (k+1)}{\tau}$ decrease it. For $\rho_\tau<0$, the argument is the same, but the eigenvalues that are dampened are those with phases $\frac{2\pi k}{\tau}$ and the enhanced ones correspond to phases $\frac{2\pi (k+1)}{\tau}$.

\begin{figure}[h!]
	\centering
	\includegraphics[width=1\linewidth,trim={0cm 1cm 0cm 1cm},clip]{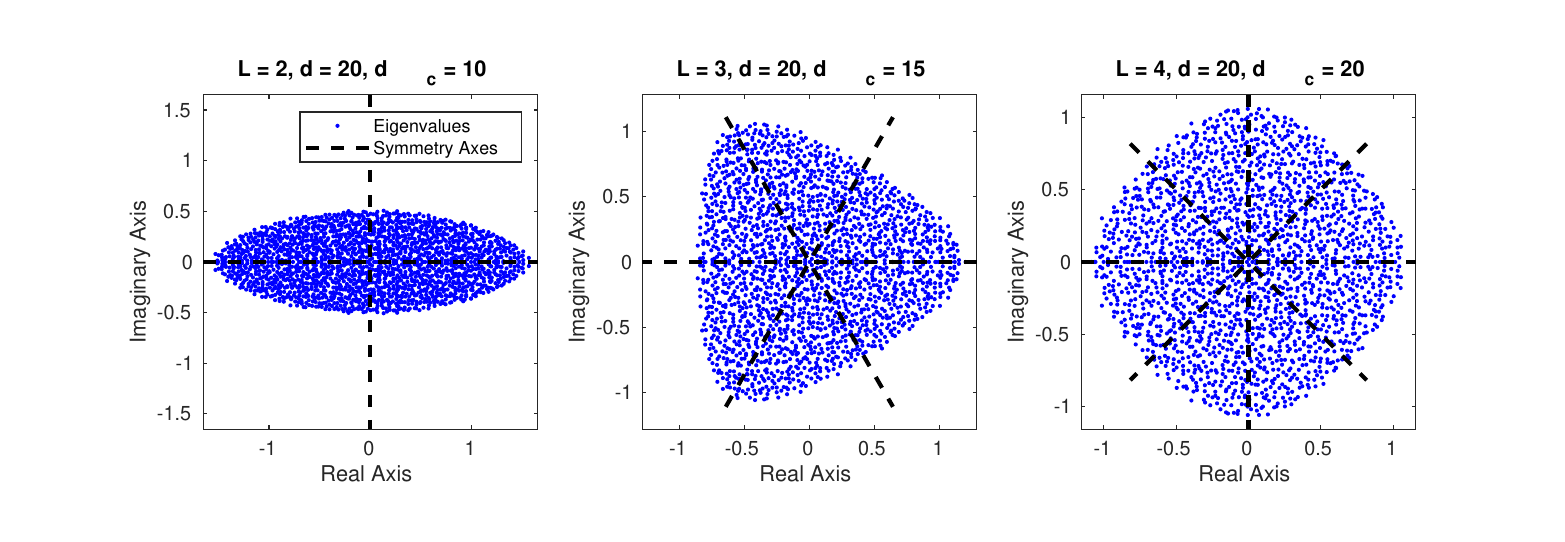}
	\caption{Eigenvalues of networks with short cycles: Each plot shows the spectrum of a single random graph with cycles of length $\tau$ being overrepresented. Each blue point corresponds to one eigenvalue in the complex plane and the black lines correspond to the axis of symmetry that we investigated. The graphs have $5000$ nodes, each with average degree of $20$, the connections are assigned randomly and the weights can be either $\frac{1}{\sqrt{d}}$ or $-\frac{1}{\sqrt{d}}$ with equal probability. }
	\label{fig:kellipseEigenvalues}   
\end{figure}

We computed the eigenvalue distributions of large graphs with overrepresented circular motifs where degrees and weights are homogeneous. The resultant eigenvalue distributions shown in Fig.\ref{fig:kellipseEigenvalues} agree with our previous discussion. In fact, their support is given by the Hypotrochoidic Law of Random Matrices \cite{aceituno2019universal}, which states that the support of the eigenvalues of such a matrix with a variance across columns normalized to 1 has the boundary
\begin{equation}
	z(\varphi) = e^{-i\varphi} + \rho_\tau e^{i(\tau -1)\varphi},
\end{equation}
where $\varphi$ is an angle. 

\section{Circulant directed networks}
\begin{figure}\label{fig:circulant}
	\centering
	\includegraphics[width=0.5\linewidth, trim={1cm 1cm 1cm 1.3cm},clip]{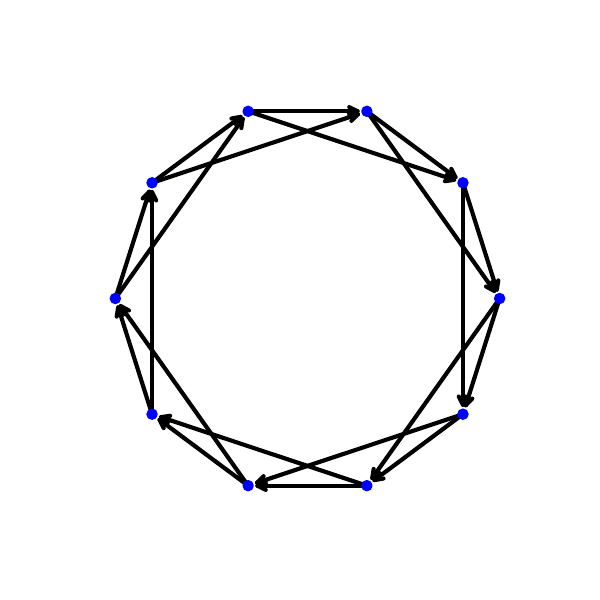}
	\caption{Circulant graph: Graphical representation of a directed circulant with graph ten nodes and degree two. % The edges of every node go to its two clockwise neighbors.
		}
\end{figure}

Having studied the spectra of random graphs with short cycles, we now investigate the effects of long cycles. We turn to circulant graphs, a well-known class of graphs that provided landmark results in graph theory \cite{watts1998collective}. In our case we take their directed version, where every node with index $n\in \mathbb{Z}/n\mathbb{Z}$ and degree $d$ is connected to its neighbors $\left[n+1, n+2, ..., n+d\right]$ as shown in Figure \ref{fig:circulant} for $d=2$. For simplicity, we consider random weights independently sampled from the set $\{-1, 1\}$ with equal probability. When $N\gg d$, this family of graphs contain only long cycles and is therefore well suited for our question. 

The simplest case is given when $d=1$. Since there is only a single simple cycle of length $N$, $\rho_k=0$ for all $k<N$, and $\rho_N = \prod_{n=1}^{N}w_{n, n+1}$. The solutions for those equations are the $N$th roots of $\rho_N$. For $\rho_N=1$,
\begin{equation}\label{eq:singleCycleEigs}
\lambda_n = e^{\frac{2\pi n}{N}i}\quad \forall n\in \mathbb{Z}/n\mathbb{Z},
\end{equation}
where $i$ is the imaginary unit. If $\rho_N<0$ the phase is shifted by $\frac{\pi}{N}$. 

\begin{figure}[h!]
	\centering
	\includegraphics[width=1\linewidth, trim={2.5cm 0cm 2cm 0cm},clip]{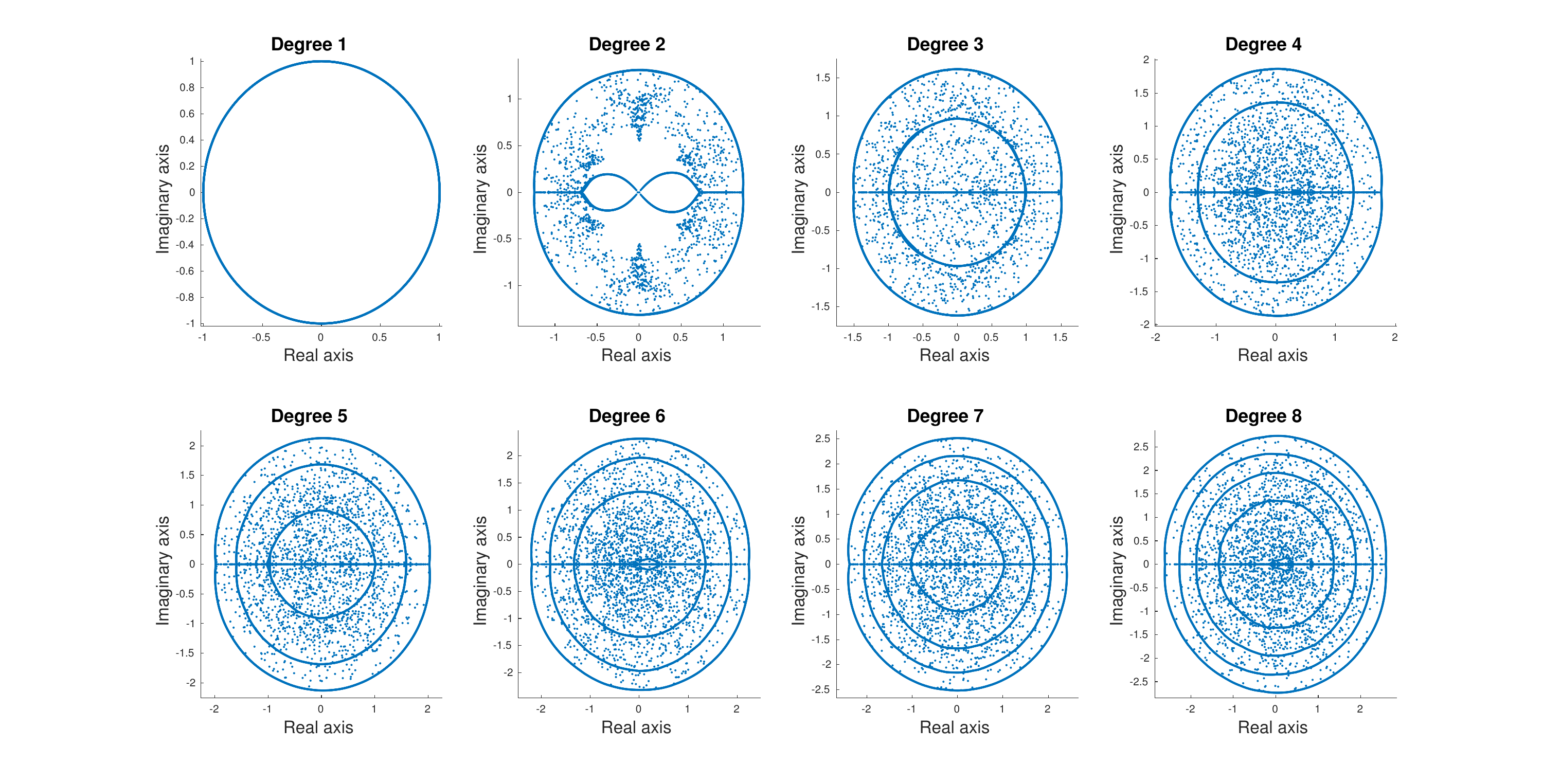}
	\caption{Eigenvalues of circulant networks: Each plot shows the spectrum of a directed circulant graph with $5000$ nodes where the weights are randomly chosen to be $-1$ or $1$ with equal probability. Each blue point corresponds to one eigenvalue in the complex plane. }
	\label{fig:circulantEigenvalues}   
\end{figure}

We computed the eigenvalues of large circulant directed networks, our results are in Fig \ref{fig:circulantEigenvalues}. For $d=1$, our prediction matches the distribution. When we increase $d$ we find a surprising result: For $d=2$, there is an outer ring and the rest of eigenvalues concentrate in a lemniscate shape, and for $d>2$ the eigenvalues are concentrated on the real line and on $\lceil\frac{d}{2}\rceil$ rings centered at the origin.

To extend Eq. \ref{eq:singleCycleEigs} to $d>1$, we start by studying the value of the largest eigenvalues by a combinatorial argument. When we take Eq. \ref{eq:cLIsSumOfEigs} in the limit where $L\rightarrow \infty$, the eigenvalues with the largest modulus dominate the rest. Then,
\begin{equation}\label{eq:lambdaToLambdaMax}
|\rho_L| = \dfrac{1}{N} |\sum_{k=1}^{N} \lambda_k^L| \sim |\lambda_\text{Max}|^L.
\end{equation}
so $|\rho_L|$ grows exponentially in $L$ at the rate $\lambda_\text{Max}$.

We can use Eq. \ref{eq:defcL} to write $\rho_L$ as a binomial distribution where cycle weights $w_c$ are independently sampled from $\{-1,1\}$ with equal probability. Then, $\rho_L$ can be approximated by a normal distribution $\mathcal{N}(0,\sigma^2_L)$ with
\begin{equation}
\sigma^2_L = \dfrac{1}{N} |C_L|,
\end{equation} 
where $|C_L|$ is the number of cycles of length $L$. To obtain this value we consider the recurrence relation between paths in the circulant graph
\begin{equation}
p_{L}(n,m) = \sum_{k=1}^{d} p_{L-1}(n,m-k),
\end{equation}
where $p_{L-1}(n,m)$ is the number of paths from node $n$ to $m$ of length $L$. When $L$ becomes large, $p_L(n,n)\sim p_L(n-k,n)$ for $k\ll N$. Therefore,
\begin{equation}
|C_L| = p_L(n,n)= \sum_{k=1}^{d} p_{L-1}(n,m-k) \sim d p_{L-1}(n,n) = d |C_{L-1}|.
\end{equation}

Given a normal distribution, the expected value of the first order-statistic has the order of magnitude of the standard deviation, so
\begin{equation}\label{eq:cLtoSqrtD}
|\rho_L| \sim \sqrt{d^{L}},
\end{equation} 
By putting together Eq. \ref{eq:lambdaToLambdaMax} and Eq. \ref{eq:cLtoSqrtD},
\begin{equation}\label{eq:specRadLimitCirculant}
|\lambda_\text{Max}| \approx \sqrt{d}.
\end{equation} 

\begin{figure}[h!]
	\centering
	\includegraphics[width=0.8\linewidth, trim={0cm 0cm 0cm 0cm},clip]{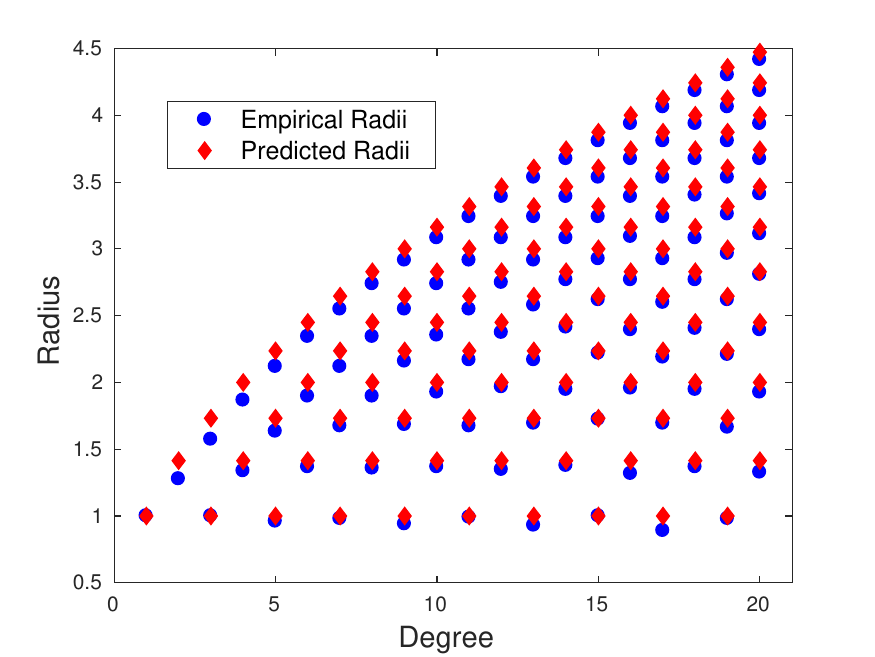}
	\caption{Eigenvalue radii for circulant random graphs: Each blue point corresponds to the radius of a circle in the distribution of eigenvalues for a circulant matrix with $5000$ nodes and signs sampled at random (see Fig. \ref{fig:circulantEigenvalues}). Every red rhombus is the radius given by Eq. \ref{eq:CirculantRadii}, and the highest ones for every degree correspond to the radius given by Eq. \ref{eq:cLtoSqrtD}.}
	\label{fig:radiiCircles}   
\end{figure}

Eq. \ref{eq:specRadLimitCirculant} gives us the approximate position of the larger ring which agrees with the observed values (see Figure \ref{fig:radiiCircles}). We also observe empirically that the radii of the rings are distributed following the rule
\begin{equation}\label{eq:CirculantRadii}
	r_d(k) \lesssim \begin{cases}
	\sqrt{2k-1} \quad &\forall \ d = 1 \mod{2}\\
	\sqrt{2k} \quad &\forall \ d = 0 \mod{2}
	\end{cases}
\end{equation}
where $r_d(k)$ is the radius of the $k$th ring in a circulant graph with degree $d$. We compare our observation with the empirical data in figure (\ref{fig:radiiCircles}), finding that the fit is quite good.

\section{Discussion}

In this paper we have argued that the structure of a graph can reveal information about its spectra by counting cycles, and shown two examples of families of matrices where some properties of their eigenvalue distribution can be studied using this relationship. Although we restricted our examples to simple graphs, we are confident that our framework can be applied to more complicated weight distributions and more complex structural features.
 
We believe that the results presented here are easily translated into complex systems setting: Cycles are easy to understand and quantify, eigenvalues are better suited to study a systems' stability and dynamics. By establishing a link between both, we provide a better insight into the relationship between both cases. Furthermore, we also showed two beautiful families of eigenvalue distributions, hopefully making random matrix theory more attractive by sheer aesthetics.

\section*{Acknowledgement}

P.V.A is supported by the Supported by BMBF and Max Planck Society

\bibliographystyle{unsrt}
\bibliography{cyclesAndSpectraBib.bib}

\end{document}